\documentclass[prepint,review,12pt,3p]{elsarticle}

\usepackage{stmaryrd}
\usepackage{amsmath}
\usepackage{amsfonts}
\usepackage{amsthm}

\journal{Journal of Computational and Applied Mathematics}

\newtheorem{lem}{Lemma}
\newtheorem{theo}{Theorem}
\newtheorem{prop}{Proposition}

\begin{document}

\begin{frontmatter}

\title{A spectral method based on $(0,2)$ Jacobi polynomials. Application to Poisson problems in a sphere.}
\author{Cornou Jean-Louis}
\ead{jean-louis.cornou@obspm.fr}
\author{Bonazzola Silvano}
\ead{silvano.bonazzola@obspm.fr}

\address{Laboratoire Univers et Th\'eories, Observatoire de Paris, CNRS, Universit\'e Paris Diderot, 5 place Jules Janssen, F-92190, Meudon, France.}

\begin{abstract}
A new spectral method is built resorting to $(0,2)$ Jacobi polynomials. We describe the origin and the properties of these polynomials. This choice of polynomials is motivated by their orthogonality properties with the respect to the weight $r^2$ used in spherical geometry. New results about Jacobi-Gauss-Lobatto quadratures are proven, leading to a discrete Jacobi transform. Numerical tests for Poisson problems in a sphere are presented using the C++ library \textsc{lorene}.
\end{abstract}

\begin{keyword}
Spectral methods \sep Jacobi polynomials \sep Numerical Relativity
\end{keyword}

\end{frontmatter}

\section{Introduction}
\label{s:intro}

Amongst the various techniques used to discretize partial differential equations, spectral methods, introduced by D. Gottlieb and S. Orszag in the seventies~\cite{gottlieborszag}, offer high accuracy at a low computational cost. Their principle is to approximate solutions of PDEs by truncated Fourier series or high degree polynomials. These methods have an infinite order of convergence, because the error between the exact and discrete solutions is only limited by the regularity of the exact solution. 

Spectral methods resort to Fourier series in the case of periodic boundary conditions; whereas polynomial approximations use orthogonal polynomials (such as Chebyshev or Legendre polynomials). These methods also employ quadrature formulas to compute integrals in variationnal formulations of PDEs. 

In the present work, we explore a new family of orthogonal polynomials : the $(0,2)$ Jacobi polynomials. After a brief summary about spectral methods in Sec.~\ref{s:spectral}, we detail several general properties of the Jacobi polynomials, and present the necessary tools to build a spectral method based on the $(0,2)$ Jacobi polynomials in Sec.~\ref{s:jacobi}. Finally, we show some applications to Poisson problems including numerical tests in Sec.~\ref{s:test}. Concluding remarks are given in Sec.~\ref{s:concl}, and important technical proofs can be found in the appendix (Sec.~\ref{s:annex}).

\section{Spectral methods}
\label{s:spectral}

Both numerical analysis and implementation of spectral methods are based on orthogonal polynomials, whose major properties are hereby recalled. We shall then emphasize the importance of Sturm-Liouville problems in the case of Legendre polynomials, and mention results about polynomial approximation and interpolation. Several textbooks cover this domain, for instance~\cite{szego},~\cite{boyd}, or~\cite{qss}, and provide comprehensive proofs of the following results.

\subsection{Orthogonal polynomials}

Let $\Lambda$ be the interval $(-1,1)$ and $w$ denote a weight on $\Lambda$, i.e, $\forall n \: \int_{-1}^1 w(x) x^n dx < \infty$ and $w>0$ on $\Lambda$. We define 
\begin{equation}
\label{e:spaceL2}
L_w^2(\Lambda) = \left\{ 
			u : \Lambda \longrightarrow \mathbb{R} \; \lvert \;
			\int_{-1}^1 u^2(x) w(x) dx < \infty 
		\right\},
\end{equation} which is a Hilbert space for the scalar product 

\begin{equation}
\label{e:scalprod}
(f \lvert g ) = \int_{-1}^1 f(x)g(x)w(x)dx  .
\end{equation}

One can construct a basis of monic orthogonal polynomials by using a Gram-Schmidt orthogonalisation process on the basis $x^n$, $n \geq 0$ :  $P_0$ is fixed to $1$, then, assuming the $P_m$, $0 \leq m \leq n-1$ are known, $P_n$ is chosen by 

\begin{equation}
\label{e:gramprocess}
P_n(x) = x^n - \sum_{m=0}^{n-1} 
	\frac{\int_{-1}^1 y^mP_m(y)w(y)dy}{\|P_m\|_{L_w^2(\Lambda)}^2} P_m(x) .
\end{equation}

\smallskip

We recall the following property of orthogonal polynomials:

\begin{prop}
For all positive integer $n$, the roots of $P_n$ are real, distinct and strictly bounded by $-1$ and $1$.
\end{prop}
 
In particular, those polynomials do not vanish at $1$, we can thus define any family of orthogonal polynomials by imposing their value at $1$. Besides, the monic orthogonal polynomials satisfy the following induction formula.

\begin{prop} 
For all integer $n \geq 2$, 
\begin{equation} \label{e:inductionorthogonal}
 P_{n}(x) = \left( x - \frac{(xP_{n-1}|P_{n-1})}{\|P_{n-1}\|^2} \right) P_{n-1}(x) -\frac{\|P_{n-1}\|^2}{\|P_{n-2}\|^2} P_{n-2}(x) .
\end{equation}
\end{prop}

\smallskip

Finally, a wide class of orthogonal polynomials belongs to singular Sturm-Liouville solutions, on which spectral discretizations rely. We present the Legendre case ($w=1$ and $L_n(1)=1$), since a general description would be too long.
 
\begin{theo}
For all integer $n \geq 0$, $L_n$ satisfies the following differential equation

\begin{equation}\label{eq-dif-legendre}
\frac{d}{dx} \left( (1-x^2)L_n'(x) \right) + n(n+1) L_n(x) = 0 .
\end{equation}

\end{theo}

The $A$ operator, defined by
\begin{equation}
\label{ope-sturm-liouville}
A \varphi = -\frac{d}{dx} \left( (1-x^2)\varphi'(x) \right)
\end{equation}
is self-adjoint in $L^2(\Lambda)$ by integration by parts. Besides, it is positive and of singular Sturm-Liouville type. The equation (\ref{eq-dif-legendre}) means that Legendre polynomials are eigenvectors of this operator, which is the origin of the  ``spectral" adjective of the numerical methods hereby discussed.

A consequence of (\ref{eq-dif-legendre}) (through integration by parts) is that, for all positive integers $m$ and $n$ :
\begin{equation}\label{norme-derivee-legendre}
\int_{-1}^1 L_m'(x)L_n'(x)(1-x^2)dx = n(n+1) \int_{-1}^1 L_m(x)L_n(x)dx.
\end{equation}
This means that the $L_n'$, $n \geq 1$ form a basis of orthogonal polynomials in $L_{1-x^2}^2(\Lambda)$. 

\subsection{Polynomial approximation error and Sturm-Liouville problem}

This section illustrates why spectral methods are so accurate. The Sturm-Liouville operator $A$ enables us to put some upper bounds on the distance between functions with a given regularity and a polynomial space, still in the case of Legendre polynomials. This distance is computed, by mean of orthogonal projectors. It has been first estimated in~\cite{cq82} and~\cite{maday90}.

Let $N$ be a positive integer, $\mathbb{P}_N(\Lambda)$ stands for the space of polynomials with degree less than $N$, and the orthogonal projector from $L^2(\Lambda)$ to $\mathbb{P}_N(\Lambda)$ is denoted by $\pi_N$. For any integer $m \geq 0$, $H^m(\Lambda)$ stands for the Sobolev space of order $m$ on $\Lambda$.

\begin{theo}\label{t:polapprox}
For all integer $m \geq 0$, there exists a positive constant $c$, only dependant of $m$ such that, for all $\varphi$ in $H^m(\Lambda)$, 
\begin{equation}
\| \varphi - \pi_N \varphi \|_{L^2(\Lambda)} \leq c N^{-m} \| \varphi \|_{H^m(\Lambda)}.
\end{equation}
\end{theo}

The key of the proof of this theorem is the following one. We denote $\varphi_n$ the $n$-th coefficient of $\varphi$ in its expansion onto the Legendre basis. Then, 
\begin{equation}
 \varphi_n = \frac{1}{\|L_n\|_{L^2(\Lambda)}^2} (\varphi \lvert L_n) = 
\frac{1}{\| L_n \|_{L^2(\Lambda)}^2} \frac{1}{n(n+1)} (\varphi \lvert AL_n).
\end{equation}
Since $A$ is self-adjoint in $L^2(\Lambda)$, 
\begin{equation} \varphi_n = \frac{1}{\| L_n \|_{L^2(\Lambda)}^2} \frac{1}{n(n+1)} (A \varphi \lvert L_n).
\end{equation}
By iterating $k$ times this result, one can deduce that
\begin{equation} \varphi_n = \frac{1}{\| L_n \|_{L^2(\Lambda)}^2} \frac{1}{\left( n(n+1) \right)^k} (A^k \varphi \lvert L_n).
\end{equation}
Since $(A^k \varphi \lvert L_n) / \| L_n \|_{L^2(\Lambda)}^2$ is the $n$-th coefficient of $A^k \varphi$ in its decomposition over the Legendre basis, a continuity result on the operator $A$ will allow to conclude. 

\subsection{Polynomial interpolation error}
Thus we have seen that any function, given its regularity, can be approximated by polynomials, with an error decreasing as a power law of the degree of the polynomial. The power depends only on the chosen norm, and the regularity of the function. However, this result has few direct numerical applications. Indeed, one has to compute integrals to obtain the polynomial approximating the function, which is highly time-expensive. Thus, one shall approximate these integrals with quadrature formulas, and replace the orthogonal projectors by interpolation operators $i_N$ in the nodes of these quadrature formulas. First, we recall the definition of these operators.
\begin{prop}\label{p:definterp}
Given $N+1$ distinct points $(x_i)$ in $\Lambda$, and a continuous function $\varphi$ on $\Lambda$, there exists a unique polynomial $i_N \varphi$ in $\mathbb{P}_N$ such that
\begin{equation}
\forall i, \; \; \; i_N\varphi (x_i) = \varphi(x_i)
\end{equation}
\end{prop}
This greatly simplifies the set up of the method, and doesn't reduce the precision according to the following theorem, first demonstrated in~\cite{maday91}. 
\begin{theo}\label{t:polinterp}
For all integer $m \geq 1$, there exists a positive constant $c$, depending only in $m$, such that, for all function $\varphi$ in $H^m(\Lambda)$, we have
\begin{equation} \| \varphi - i_N \varphi \|_{L^2(\Lambda)} \leq c N^{-m} \| \varphi \|_{H^m(\Lambda)} . \end{equation}

\end{theo}

The latter result shows spectral methods have an infinite order of accuracy. A $\mathcal{C}^{\infty}$ function is interpolated faster than any power of the discretization parameter. Moreover, it can be shown that an exponential convergence of the interpolant is achieved if the function is analytic. 

\section{Jacobi polynomials}
\label{s:jacobi}

In this section, several results about Jacobi polynomials are detailed. First, we present how Jacobi polynomials arise naturally from Sturm Liouville singular problems. Then, after presenting some basic properties, we prove a new result about Jacobi-Gauss-Lobatto quadratures, enabling us to build a discrete Jacobi transform. Finally, we show some differentiation and integration results about $(0,2)$ Jacobi polynomials.

\subsection{Introducing Jacobi polynomials}

We saw previously that the Liouville operator was crucial to obtain an efficient polynomial approximation. More generally, a Sturm-Liouville problem consists in looking for solutions $(u,\lambda)$ of

\begin{equation} \left\{ \begin{array}{c c}  -(pu')'+qu = \lambda w u  \text{  on the interval  } (-1,1) 
			\\  \text{Suitable boundary conditions for } \, u 
		\end{array} \right. \end{equation}

The coefficients $p$, $q$ and $w$ are three given, real-valued functions such that $p$ is continuously differentiable, strictly positive in $(-1,1)$ and continuous at $x=\pm1$; $q$ is continuous, non-negative and bounded in $(-1,1)$; the weight $w$ is continuous and non-negative in $(-1,1)$. One must notice that every Sturm-Liouville problem does not ensure the convergence of the associated spectral method. Let us consider the following example ;

\begin{equation} \left\{ \begin{array}{l l}  u''+\lambda u = 0 
		\\ u'(-1)=u'(1)=0 \end{array} \right. \end{equation}

Its eigenvalues are $\lambda_k=(\pi k)^2/2$, with eigenfuctions $\phi_k(x)=\cos((\pi/2) k (x+1))$. Thus, a function will be correctly approximated by the cosine series, iff all its even derivatives vanish in $-1$ and in $1$. This is due to the non-vanishing of $p$ in the extremities of the interval.This Sturm-Liouville problem is \emph{regular}. Inversely, if $p$  vanishes in $-1$ and in $1$, the problem is said to be \emph{singular}. In the case of singular problems, efficient polynomial approximation is ensured. We send back to~\cite{chqz1} for a more complete discussion.

Thus, we are led to consider polynomial solutions to singular Sturm-Liouville problems in order to build an efficient spectral method. Let us note $(\phi_k)$ a family of polynomial solutions of degree $k$ of a singular Sturm-Liouville problem with eigenvalues $\lambda_k$, then one can see that $q/w$ is a polynomial of degree zero, $p/w$ has degree $2$, and $p'/w$ has degree $1$, according to the following identity:

\begin{equation} \forall k \; \; \; \phi_k = -\frac{p}{\lambda_k w} \phi_k'' - \frac{p'}{\lambda_k w} \phi_k' + \frac{q}{\lambda_k w} \phi_k. \end{equation}

Since $p$ vanishes in $-1$ and in $1$, $p/w= c_0(1-x)(1+x)$. Let us note $ax+b=p'/w$, then $p'/p = (ax+b)/[(1-x)(1+x)]=a_1/(1-x)+b_1/(1+x)$. After integration, one has $p=c_1(1-x)^{\alpha+1}(1+x)^{\beta+1}$ with $c_1$, $\alpha$ and $\beta$ constants. As a result, $w=c_2(1-x)^{\alpha}(1+x)^{\beta}$. An integrability condition on $w$ imposes $\alpha > -1$ and $\beta > -1$.

We then define the Jacobi polynomials $J_n^{(\alpha,\beta)}$ of index $(\alpha,\beta)$ as the orthogonal polynomials for the weight $w(x)=(1-x)^{\alpha}(1+x)^{\beta}$, normalized by \begin{equation} \label{e:normalisationjacobi}
J_n^{(\alpha,\beta)}(1)=\frac{\Gamma(n+1+\alpha)}{\Gamma(1+\alpha)\Gamma(n+1)},\end{equation} where $\Gamma$ is the Euler gamma function. 

These polynomials are solutions of the previous singular Sturm-Liouville problem. Indeed, $\frac{1}{w} \frac{d}{dx} \left(  (1-x^2) w \frac{d}{dx}J_n^{(\alpha,\beta)}  \right)$ is a polynomial with degree $\leq n$, and verifies that for all polynomial $\varphi$ with degree $\leq n-1$

\begin{equation} \int_{-1}^1 \frac{1}{w}\frac{d}{dx} \left(  (1-x^2) w \frac{d}{dx} J_n^{(\alpha,\beta)} \right) \varphi w = \int_{-1}^1 J_n^{(\alpha,\beta)}  \frac{1}{w} \frac{d}{dx} \left(  (1-x^2) w \varphi' \right) w = 0 .\end{equation}
Therefore there exists $\lambda_n$ such that 
\begin{equation} \label{eqdifjacobi} \frac{1}{w} \frac{d}{dx} \left(  (1-x^2) w \frac{d}{dx} J_n^{(\alpha,\beta)}  \right) = - \lambda_n J_n^{(\alpha,\beta)} .\end{equation}
Jacobi polynomials are thus the eigenvectors of the previous Sturm-Liouville problem. Identifying the highest degree coefficients, one obtains $\lambda_n = n(n+\alpha+\beta+1)$. Thus, this configuration enables to generalize Theorem \ref{t:polapprox}. 

Standard orthogonal polynomials are special cases of Jacobi polynomials. Legendre polynomials are Jacobi polynomials with index $\alpha=\beta=0$. If we denote $T_n(x) = \cos(n\arccos(x))$ the $n$-th Chebyshev polynomial, one can verify that \begin{equation}\label{e:tcheb} T_n=J_n^{(-1/2,-1/2)}\Gamma(n+1)\Gamma(1/2)/\Gamma(n+1/2).\end{equation} 
Finally, $(0,1)$ Jacobi polynomials have been studied in~\cite{bernardidaugemaday}.

\subsection{Properties of Jacobi polynomials}

We present here basic results about Jacobi polynomials. It can be proven they are special cases of hypergeometric functions (see~\cite{handbook}), which provides the following results :
\begin{itemize}
\item Analytical expression
\begin{equation} J_k^{(\alpha,\beta)}(x)= \frac{1}{2^k} \sum_{l=0}^k \left(\!\!\!
  \begin{array}{c}
    k+\alpha \\
    l
  \end{array}
  \!\!\!\right) \left(\!\!\!
  \begin{array}{c}
    k+\beta \\
    k-l
  \end{array}
  \!\!\!\right)
(x-1)^l (x+1)^{k-l} .\end{equation}

\item Induction formula
\begin{multline} \label{inducjacobi}
2(k+1)(k+\alpha+\beta+1)(2k+\alpha+\beta)J_{k+1}^{(\alpha,\beta)}(x)
\\
= [(2k+\alpha+\beta+1)(\alpha^2 - \beta^2) + x\Gamma(2k+\alpha+\beta+3)/\Gamma(2k+\alpha+\beta) ]J_{k}^{(\alpha,\beta)}(x)
\\
- 2(k+\alpha)(k+\beta)(2k+\alpha+\beta+2) J_{k-1}^{(\alpha,\beta)}(x).
\end{multline}

\item Highest degree coefficient
 
\begin{equation} k_n = \frac{1}{2^n}\frac{\Gamma(2n+\alpha+\beta+1)}{\Gamma(n+1)\Gamma(n+\alpha+\beta+1)} .\end{equation}

\bigskip

\item Norm

\begin{equation} \label{e:jacobinorm} \| J_n^{(\alpha,\beta)}\|_{L_w^2(\Lambda)}^2 = \frac{2^{\alpha+\beta+1} \Gamma(n+\alpha+1) \Gamma(n+\beta+1)}{(2n+\alpha+\beta+1)\Gamma(n+1)\Gamma(n+\alpha+\beta+1)}. \end{equation}

\bigskip

\item Links between the various Jacobi polynomials families.
\begin{equation} \label{linkjacobi} \begin{array}{l}
(n+\alpha/2 + \beta/2 +1)(1-x)J_n^{(\alpha+1,\beta)} = (n+\alpha+1)J_n^{(\alpha,\beta)} - (n+1)J_{n+1}^{(\alpha,\beta)}.
\\
(n+\alpha/2 + \beta/2 +1)(1+x)J_n^{(\alpha,\beta+1)} = (n+\beta+1)J_n^{(\alpha,\beta)} + (n+1)J_{n+1}^{(\alpha,\beta)}.
\\
(2n+\alpha+\beta)J_n^{(\alpha-1,\beta)} = (n+\alpha+\beta)J_n^{(\alpha,\beta)}-(n+\beta)J_{n-1}^{(\alpha,\beta)}.
\\
(2n+\alpha+\beta)J_n^{(\alpha,\beta-1)} = (n+\alpha+\beta)J_n^{(\alpha,\beta)}+(n+\alpha)J_{n-1}^{(\alpha,\beta)}.
\\
J_n^{(\alpha,\beta)}(-x) = (-1)^n J_n^{(\beta,\alpha)}(x).
\end{array} \end{equation}

\end{itemize}
The latter formulae show, for example, that $(0,2)$ Jacobi polynomials are linked to the Legendre polynomials via the following equation :
\begin{equation} \label{e:linkjacobilegendre}
(n+3/2)(1+x)^2 J_n^{(0,2)} = (n+2)L_n + (2n+3)L_{n+1} + (n+1)L_{n+2}
\end{equation}
Finally, as we saw it previously for Legendre polynomials, derivatives of Jacobi polynomials are orthogonal for the weight $w(x)(1-x^2)$, which can also be shown from the identity

\begin{equation} \label{derjacobi} \frac{d}{dx} \left( J_n^{(\alpha,\beta)} \right) = \frac{1}{2}(n+\alpha+\beta+1) J_{n-1}^{(\alpha+1,\beta+1)} .\end{equation}

We show in a very illustrative manner the plot of $J_{n}^{(0,2)}$, $2 \leq n \leq 5$, in Fig.1.

\begin{figure}[t]
\centerline{\includegraphics[width=\columnwidth]{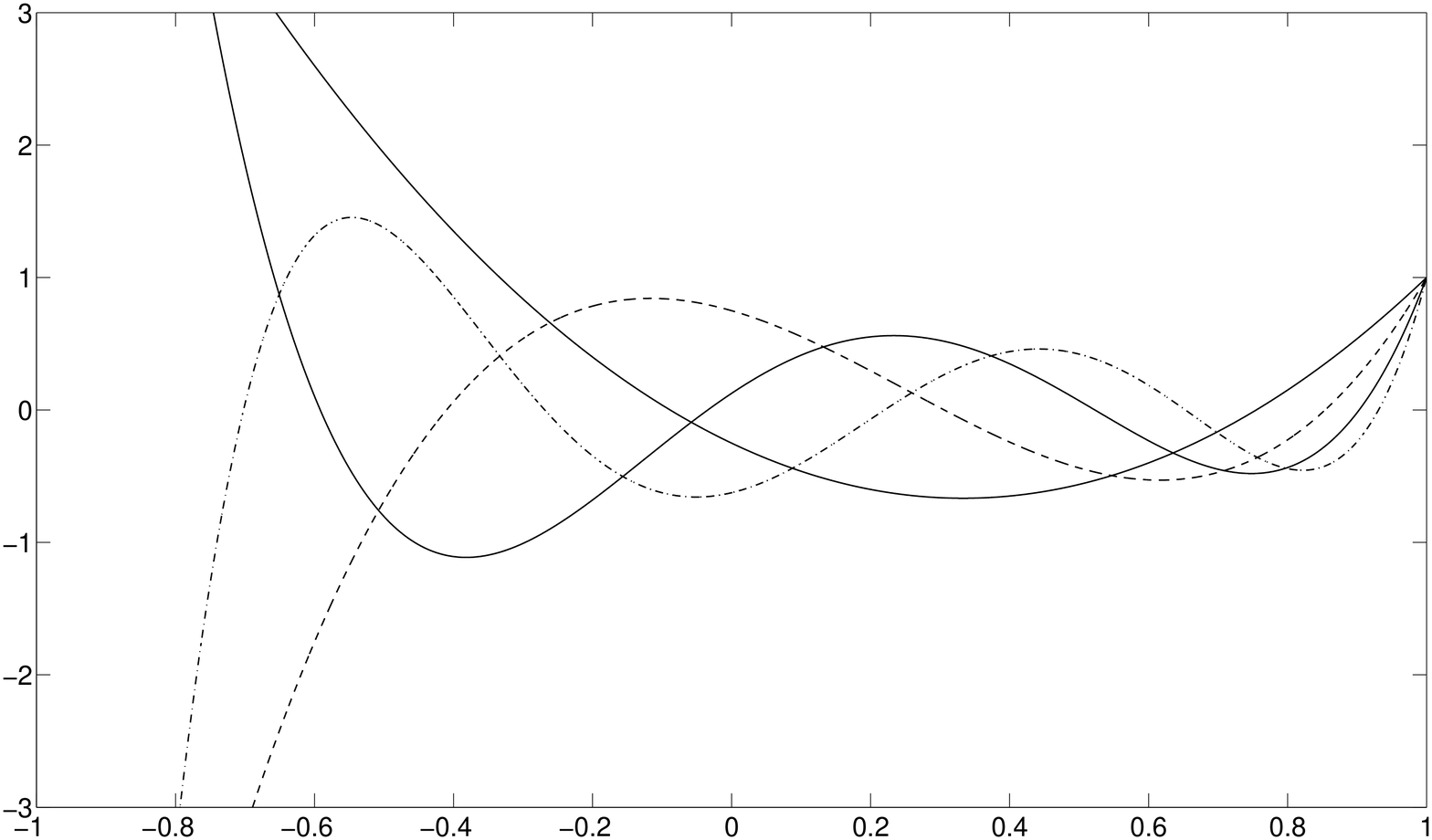}}
\caption[]{\label{f:jacobiplot} Plot of $J_{n}^{(0,2)}$, $2 \leq n \leq 5$. Notice that $J_n^{(0,2)}(-1) = (-1)^n(n+1)(n+2)/2$}
\end{figure}

\subsection{Gauss-Lobatto methods and Discrete Jacobi Transform}

The extrema of orthogonal polynomials are useful in order to construct high precision  numerical quadrature formulas, i.e, which are exact within a space of high degree polynomials : the so-called Gauss-Lobatto formulas. In this section, we recall the definition of these quadrature formulas in proposition \ref{p:defquad}, then compute some useful formulas to build numerical spectral methods.

\begin{prop} \label{p:defquad} Let $N$ be a positive integer. Let us denote $x_0=-1$ and $x_N=1$. There exists a unique set of $N-1$ points $x_i$ in $\Lambda$, $1 \leq i \leq N-1$, and a unique set of $N+1$ real $\rho_i$, $0 \leq i \leq N$, such that we have the following identity for all $\Phi$ in $\mathbb{P}_{2N-1}(\Lambda)$ 
\begin{equation} \label{quadformula} \int_{-1}^1 \Phi(x) w(x) dx = \sum_{j=0}^N \Phi(x_j) \rho_j .\end{equation}
The nodes $x_i$, $1 \leq i \leq N-1$, are the zeros of $\frac{d}{dx}J_N^{(\alpha,\beta)}$. Moreover, the $\rho_j$, $0 \leq j \leq N$ are positive.
\end{prop}

The quadrature formula
\begin{equation} \int_{-1}^1 \Phi(x) w(x) dx \simeq \sum_{j=0}^N \Phi(x_j) \rho_j \end{equation}
is called \emph{Jacobi type Gauss-Lobatto quadrature formula} with $N+1$ points. The expression of the coefficients (also called weights) $\rho_j$, $0 \leq j \leq N$ is given by the following proposition.

\begin{prop}\label{p:JGLW}

\begin{equation} \rho_0 = (\beta+1) \frac{2^{\alpha+\beta+1}}{N(N+\alpha+\beta+1)} \frac{\Gamma(N+1+\alpha)\Gamma(N+1+\beta)}{\Gamma(N+1)\Gamma(N+1+\alpha+\beta)} \frac{1}{(J_N^{(\alpha,\beta)}(x_0))^2}. \end{equation}

\begin{equation}  \forall i \in [1,N-1], \: \rho_i = \frac{2^{\alpha+\beta+1}}{N(N+\alpha+\beta+1)} \frac{\Gamma(N+1+\alpha)\Gamma(N+1+\beta)}{\Gamma(N+1)\Gamma(N+1+\alpha+\beta)} \frac{1}{(J_N^{(\alpha,\beta)}(x_i))^2}.  \end{equation} 

\begin{equation} \rho_N = (\alpha+1) \frac{2^{\alpha+\beta+1}}{N(N+\alpha+\beta+1)} \frac{\Gamma(N+1+\alpha)\Gamma(N+1+\beta)}{\Gamma(N+1)\Gamma(N+1+\alpha+\beta)} \frac{1}{(J_N^{(\alpha,\beta)}(x_n))^2}.\end{equation}
\end{prop}

\medskip

Some commentaries arise naturally at this point. As far as the authors know, this expression of the Jacobi Gauss-Lobatto weights has not been given previously in the literature. Szego, in \cite{szego}, gives a very similar expression, but only for Gauss quadrature (i.e, the nodes of the quadrature are the zeros of $J_N$). One can also check the validity of this formula for known values of $\alpha$ and $\beta$. For Legendre polynomials (i.e, $\alpha=\beta=0$), one can check that $\forall i \in [0, N] \, \rho_i = 2/[N(N+1)L_N^2(x_i)]$. As far as $(0,1)$ Jacobi polynomials (denoted $M_n$ hereafter) are concerned, our formulas coincide with those found in \cite{bernardidaugemaday}, namely $\rho_0 = 8/[N(N+2)M_N^2(x_0)]$ and $\forall i \in [1,N] \, \rho_i = 4/[N(N+2)M_N^2(x_i)] $. For Chebyshev polynomials ($\alpha = \beta = -1/2 $), this expression is also valid. Indeed, 

\begin{equation} \forall i \in [1,N-1] \: \: \: \: \rho_i = \frac{1}{N^2}\frac{\Gamma(N+1/2)^2}{\Gamma(N+1)\Gamma(N)}\frac{1}{ ( J_N^{(-1/2,-1/2)}(x_i) )^2 }. \end{equation}
Then equation (\ref{e:tcheb}) allows us to write $\rho_i = \Gamma(1/2)^2/[NT_N^2(x_i)]$. However, in the Chebyshev case, $x_i = - \cos(\pi i /N)$, so that $T_N(x_i) = (-1)^{N+i}$. 
Finally, $\forall i \in [1, N-1] \, \rho_i = \pi /N $, and $\rho_0 = \rho_N = \pi / 2N$. For the numerical applications we have in mind, $\alpha=0$ and $\beta=2$, so that $\rho_0 = 24/[N(N+3)(J_N^{(0,2)}(x_0))^2]$ and $\forall i \in [1, N] \, \rho_i = 8/[N(N+3)(J_N^{(0,2)}(x_i))^2]$.

Now, it is necessary to know how to compute the nodes of the Gauss-Lobatto formulas. Looking through the zeros of $J_N'$ via a secant method is simple, but zeros tend to accumulate near the boundaries of the interval. With no preliminary knowledge of the distribution of those knots, it is preferable to use a eigenvalue location method which is enabled by the following proposition.

\begin{prop}\label{p:matrixnodes}
The $x_j$, $1 \leq j \leq N-1$, are the eigenvalues of the tridiagonal symmetric matrix

\begin{equation} \label{e:matrixnodes} \left( \begin{array}{c c c c c}
\delta_1 & \gamma_1 & 0        & 0            & 0            \\
\gamma_1 & \delta_2 & \ddots   & 0            & 0            \\
0        & \ddots   & \ddots   & \ddots       & 0            \\
0        & 0        & \ddots   & \delta_{N-2} & \gamma_{N-2} \\
0        & 0        & 0        & \gamma_{N-2} & \delta_{N-1}
 \end{array} \right), \end{equation}
with
\begin{equation} \delta_n = -\frac{(\alpha-\beta)(\alpha+\beta+2)}{(2n+\alpha+\beta)(2n+\alpha+\beta+2)}, \quad 1 \leq n \leq N-1, \end{equation}
\begin{equation} \gamma_n = \frac{2}{2n+\alpha+\beta+2} \sqrt{\frac{n(n+\alpha+1)(n+\beta+1)(n+\alpha+\beta+2)}{(2n+\alpha+\beta+1)(2n+\alpha+\beta+3)}} , \quad 1 \leq n \leq N-2. \end{equation}
\end{prop}

The form of the matrix allows for a fast and robust computation of the eigenvalues, e.g, by a Givens-Householder algorithm. In the case of $\alpha =0$ and $\beta= 2$, 
\begin{equation} \delta_n = \frac{2}{(n+1)(n+2)}, \;\;\; \text{and} \;\;\; \gamma_n = \frac{1}{n+2}\sqrt{\frac{n(n+1)(n+3)(n+4)}{(2n+3)(2n+5)}}. \end{equation}

\medskip

Finally, these quadrature formulas, enable us to build a discrete Jacobi transform, i.e, if one knows the values of a function on the Gauss-Lobatto nodes, on can compute the coefficients of its polynomial interpolant in the Gauss-Lobatto nodes.

\begin{prop}\label{p:djt}

Let us note $i_Nf = \sum_{k=0}^N \tilde{f}_k J_k$ the interpolant of f on the Gauss-Lobatto nodes, then
\begin{multline} \forall m \leq N-1,  
\\ 
\tilde{f}_m = \frac{2m+\alpha+\beta+1}{N(N+\alpha+\beta+1)} \frac{\Gamma(m+1)\Gamma(N+1+\alpha)\Gamma(N+1+\beta)\Gamma(m+\alpha+\beta+1)}
{\Gamma(N+1)\Gamma(m+1+\alpha)\Gamma(m+1+\beta)\Gamma(N+\alpha+\beta+1)} 
 \\ 
\left\{ (1+\beta) \frac{f(x_0)J_m(x_0)}{(J_N(x_0))^2} + 
\sum_{i=1}^{N-1} \frac{f(x_i)J_m(x_i)}{(J_N(x_i))^2} +
(1+\alpha) \frac{f(x_N)J_m(x_N)}{(J_N(x_N))^2} \right\}, 
\end{multline}

\begin{equation} \tilde{f}_{N} = \frac{1}{(N+\alpha+\beta+1)} 
\left\{ (1+\beta) \frac{f(x_0)}{J_N(x_0)} + 
\sum_{i=1}^{N-1} \frac{f(x_i)}{J_N(x_i)} +
(1+\alpha) \frac{f(x_N)}{J_N(x_N)} \right\}. \end{equation}

\end{prop}

Various cases of these formulae can be checked. For Legendre polynomials, they reduce to
\begin{equation}
\tilde{f}_m = \frac{2m+1}{N(N+1)}\sum_{i=0}^N \frac{f(x_i)L_m(x_i)}{L_N^2(x_i)},
\;\;\;\;\;\;\;\;\;\;
\tilde{f}_N = \frac{1}{(N+1)}\sum_{i=0}^N \frac{f(x_i)}{L_N(x_i)}.
\end{equation}
In our case ($\alpha = 0$ and $\beta=2$), 
\begin{equation}
\tilde{f}_m = \frac{2m+1}{N(N+3)}\sum_{i=0}^N \frac{f(x_i)J_m(x_i)}{J_N^2(x_i)},
\;\;\;\;\;\;\;\;\;\;
\tilde{f}_N = \frac{1}{(N+3)}\sum_{i=0}^N \frac{f(x_i)}{J_N(x_i)}.
\end{equation}

This discrete Jacobi transform is an essential ingredient in order to build a Tau method \cite{chqz1}. Indeed, any function $f$ is represented by the coefficients $\tilde{f}_k$ from the latter proposition. Differential operators on $f$ can be translated as linear operators on the $(\tilde{f}_k)$, boundary conditions as linear constraints, etc.

Since we have detailed the properties of Jacobi-Gauss-Lobatto quadrature formulas, we may investigate the polynomial interpolation error. Theorem \ref{t:polinterp} has been proven for Legendre in~\cite{maday91}, Chebyshev in~\cite{q87}, $(0,1)$ Jacobi in~\cite{bernardidaugemaday}, and Gegenbauer polynomials ($\alpha=\beta$) in \cite{madaybernardi}. It still remains to be computed in the $(0,2)$ case. As an illustration of this section, we show the plot of $J_{12}^{(0,2)}$ and its Gauss-Lobatto nodes in Fig.2, we also include the location of Chebyshev Gauss-Lobatto nodes of the same order for comparison.

\begin{figure}[t]
\centerline{\includegraphics[width=1.2\columnwidth]{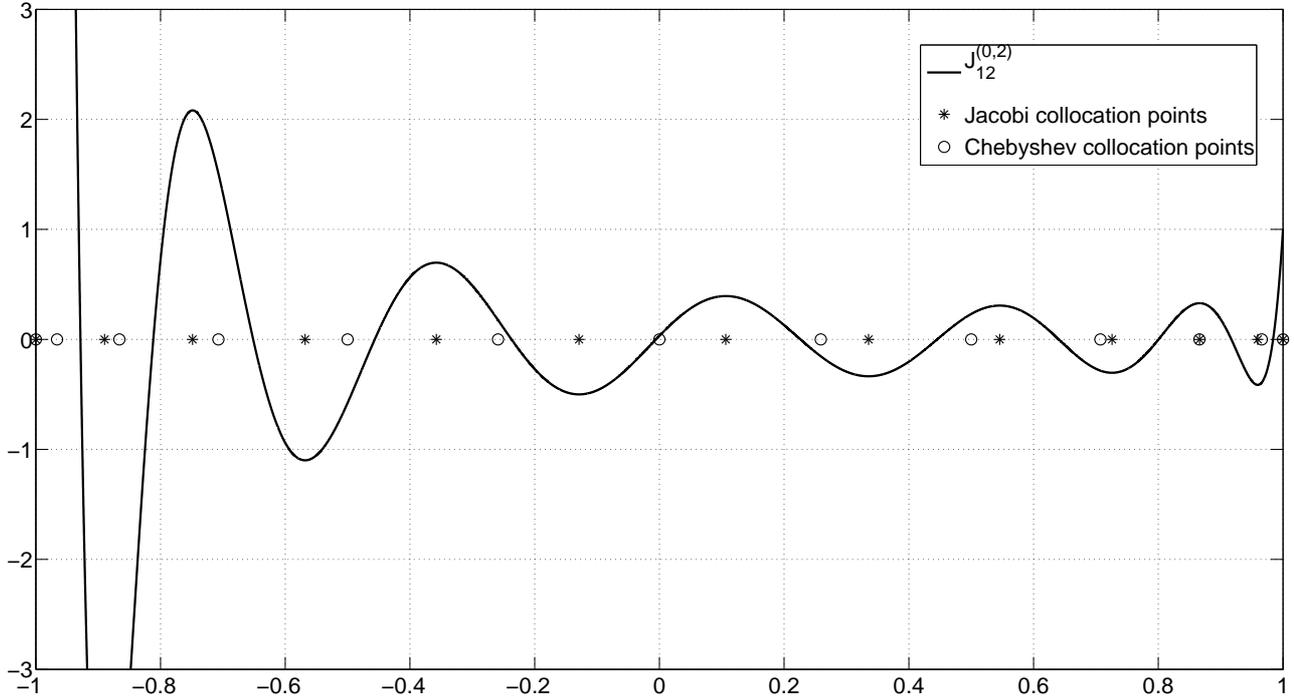}}
\caption[]{\label{f:jacobigl} Plot of $J_{12}^{(0,2)}$, location of Jacobi and Chebyshev Gauss-Lobatto nodes.}
\end{figure}

\subsection{Inverse inequalities}
The following result is enclosed here for completeness. It may not be of interest for our numerical purposes, but may be useful in order to devise some \emph{a priori} inequalities in non-linear problems. See~\cite{bernardimadayrapetti} for more details.

\begin{theo}\label{t:invineq}
There exists a real $c$ such that the following inequality is satisfied for all positive integer $N$ and all polynomial $\varphi_N$ in $\mathbb{P}_N(\Lambda)$.
\begin{equation} |\varphi_N|_{H^1_w(\Lambda)} \leq c N^2 \| \varphi_N \|_{L^2_w(\Lambda)}. \end{equation}
\end{theo}
This inequality is optimal. Indeed, there exists a constant $c'$ such that

\begin{equation} |J_N'|_{H^1_w(\Lambda)} \geq c' N^2 \| J_N' \|_{L^2_w(\Lambda)}. \end{equation}

\subsection{Physical motivation for $(0,2)$ Jacobi polynomials}
\label{motivation}

Let us have a closer look at the orthogonality property of $(0,2)$ Jacobi polynomials, which we shall denote $J_k$ hereafter.

\begin{equation} \forall n,k, \: \int_{-1}^1 J_k(x) J_n(x) (1+x)^2 dx = \delta_{k,n} \| J_k \|_w^2 .\end{equation}
By mean of the change of variable $r=(1+x)/2$, we obtain

\begin{equation} \forall n,k, \: \int_{0}^1 J_k(2r-1) J_n(2r-1) r^2 dr = \frac{1}{8} \delta_{k,n} \| J_k \|_w^2 .\end{equation}

Thus appears the crucial $r^2$ weight in spherical coordinates. The utility of such a family of polynomial appears when computing global quantities. Let us take the example of a magnetic field in a spherically symetric geometry. Its total energy is proportionnal to $ \int B^2 r^2 dr $. Denoting $b_n$ its coefficients in a decomposition on a Jacobi polynomial basis, its total energy will then be proportionnal to $\sum b_n^2 \| J_n \|_w^2 $. Using Chebyshev polynomials would have required to employ desialasing techniques to compute this integral (See \cite{chqz1}). 

Besides, when dealing with non-linear terms, one can use, for example, truncation techniques, and in our case, truncating Jacobi coefficients doesn't increase the energy of the field. For instance, if the previous magnetic field obeys an induction equation with dissipation, exact solutions have decreasing energy, and we ensure the desaliasing techniques do not increase the energy of the numerical solution.

\subsection{Operator matrices for $(0,2)$ Jacobi polynomials}

Various approaches can be examined when building a spectral method. A C++ library called \textsc{lorene} \cite{lorene} has been built and mostly uses the Lanczos-Tau method exposed in~\cite{chqz1}. Therefore, we need operator matrices such as derivation, integration, multiplication and division by $X+1$, whose expressions in the $(0,2)$ case are given in the following proposition 

\begin{prop}\label{p:opmatrices} For all positive integer $n$, the following identities are true :
\begin{itemize}
\item Derivation : 

\begin{equation} \label{e:derjacobi} J_n' = \sum_{j=0}^{n-1} \left( j+\frac{3}{2} \right) \left[ 1 - (-1)^{n-j}\frac{(j+1)(j+2)}{(n+1)(n+2)} \right] J_j .\end{equation}

\item Integration :

\begin{equation} \int_1^x J_n = \frac{n+3}{(n+2)(2n+3)}J_{n+1} - \frac{1}{(n+1)(n+2)}J_n - \frac{n}{(n+1)(2n+3)} J_{n-1} .\end{equation}

\item Division by $1+x$

\begin{equation} \label{e:divjacobi} \frac{J_n-J_n(-1)}{1+x} = \sum_{j=0}^{n-1} (-1)^{n-1-j}\frac{2j+3}{4} \left[ \frac{(n+1)(n+2)}{(j+1)(j+2)} - \frac{(j+1)(j+2)}{(n+1)(n+2)} \right] J_j .\end{equation}

\item Multiplication by $1+x$

\begin{equation} (1+x)J_n = \frac{(n+1)(n+3)}{(n+2)(2n+3)}J_{n+1} + \frac{n^2+3n+3}{(n+1)(n+2)}J_n + \frac{n(n+2)}{(n+1)(2n+3)}J_{n-1} .\end{equation}
\end{itemize}
\end{prop}

\section{Numerical tests and implementation}
\label{s:test}

Numerical tests have been performed using the C++ library \textsc{lorene}~\cite{lorene} which provides useful tools for numerical relativity, especially elliptic solvers. Poisson equations are of outmost importance when looking for solutions of Einstein equations. Indeed, various formulations of these equations lead to elliptic equations which can be solved through iterative Poisson resolutions. 

\textsc{lorene} performs 3D multi-domain resolutions with spectral accuracy thanks to Chebyshev polynomials in the radial direction, and spherical harmonics or trigonometric polynomials in the $\theta$ and $\varphi$ directions. Moreover, it uses several spherical domains, a nucleus centered around $r=0$, some shells, and an external compactified domain (using $u=1/r$ as a variable). Some applications can be found in~\cite{bgm}, \cite{gbgm}, and \cite{ncv}.

We implemented the use of the $(0,2)$ Jacobi polynomials in the nucleus, where their orthogonality properties are the most interesting, and developped a Poisson solver thanks to these polynomials.

\begin{figure}[t]
\centerline{\includegraphics[width=0.85\columnwidth]{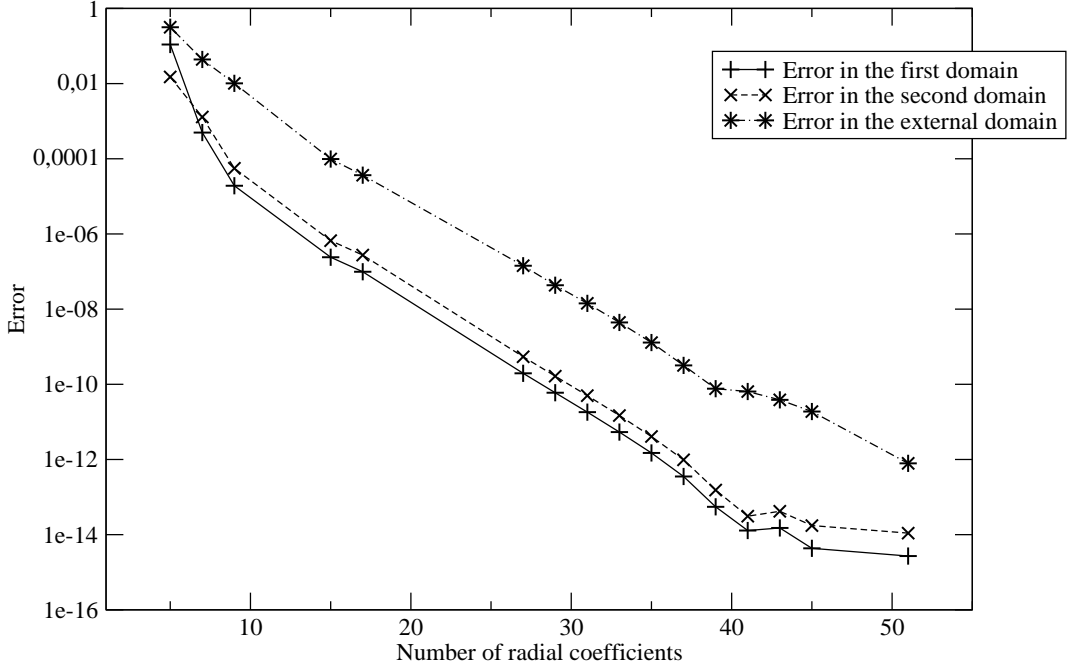}}
\caption[]{\label{f:CV_SM} Decay of the errors (maximum difference with theoretical
  solution on collocation points) in each domain. Settings are $N_{\theta}=17, N_{\varphi}=16$.}
\end{figure}

The Poisson equation reads
\begin{equation} \label{e:poisson}
\Delta f = S,
\end{equation}
where $\Delta= \partial_r^2 + \frac{2}{r} \partial_r + 1/r^2(\partial_{\theta}^2 + \frac{\cos\theta}{\sin\theta} \partial_{\theta} + \partial_{\varphi}^2)$ is the Laplace operator. Our procedure is the following one : we perform a decomposition of functions on spherical harmonics, through
\begin{equation}
f(r,\theta,\varphi)= \sum_{l,m} f_{lm}(r) Y_l^m(\theta,\varphi),
\end{equation}
where $Y_l^m$ denote the spherical harmonics. Since $\Delta Y_l^m = -l(l+1) Y_l^m$, we reduce equation (\ref{e:poisson}) to a set of 1D equations
\begin{equation} \label{e:poisson1D}
\frac{d^2 f_{lm}}{dr^2} + \frac{2}{r}\frac{d f_{lm}}{dr} - \frac{l(l+1)}{r^2}f_{lm} = S_{lm}.
\end{equation}

Then, assuming $f_{lm}(r) = \sum_{k=0}^{N_r} \tilde{f}_{klm} J_k(r)$, one can assemble the matrix of the operator appearing in (\ref{e:poisson1D}), and adding some boundary conditions, invert it to find the coefficients $\tilde{f}_{klm}$. More precisely, we decompose $\mathbb{R}^3$ in three domains, a sphere of radius $1$, a shell between radii $1$ and $2$, and an external zone $r \geq 2$. We map those three domains to $[-1,1]$ via the following equations :
\begin{equation}
\begin{array}{l l r}
\text{for} & r \leq 1, & r=(1+x)/2, \\
\text{for} & 1 \leq  r \leq 2, & r = (3+x)/2, \\
\text{for} & 2 \leq r, & r = 4/(1-x).
\end{array}
\end{equation}
Then, for each $l$, we assemble a three-block matrix containing the matrices of the following operator in the $(0,2)$ Jacobi basis for the nucleus thanks to (\ref{e:derjacobi}) and (\ref{e:divjacobi}), and in the Chebyshev basis for the two other domains 

\begin{figure}[t]
\centerline{\includegraphics[width=0.85\columnwidth]{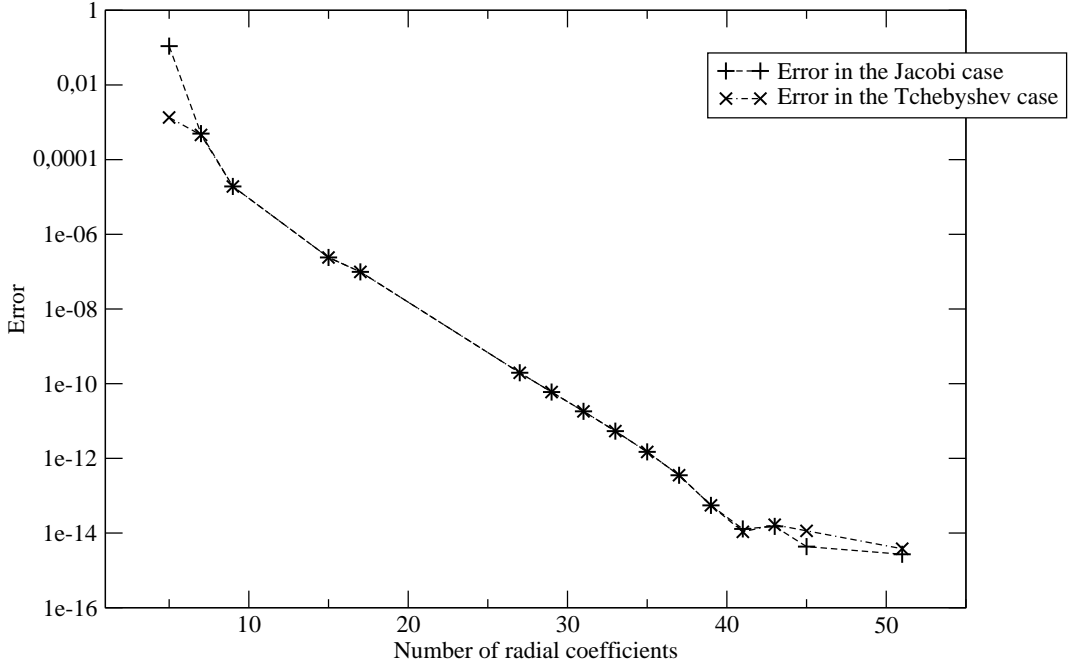}}
\caption[]{\label{f:CV_CP} Decay of the errors (maximum difference with theoretical
  solution on collocation points) in the nucleus for Jacobi and Chebyshev methods. Settings are $N_{\theta}=17, N_{\varphi}=16$.}
\end{figure}

\begin{equation}
\frac{d^2 f}{dx^2} + \frac{2}{1+x}[\frac{df}{dx} - f'(-1)] - \frac{l(l+1)}{(1+x)^2}[f - f(-1) - (1+x)f'(-1)].
\end{equation} 
In the nucleus, it corresponds to the operator in (\ref{e:poisson1D}), under the assumption that for $l=0$, $f'(-1)=0$, for $l=1$, $f(-1)=0$, and $\forall l \geq 2$, $f(-1)=f'(-1)=0$. In order to enforce those conditions when solving the problem, we replace the last lines of the first block-matrix by the values of $J_n(-1)$ or $J_n'(-1)$, and those of the source by $0$, depending of the value of $l$. Finally, we include $\mathcal{C}^1$ matching conditions between the domains in the 3-block matrix.

Several numerical tests have been passed. First, we used the following smooth source in Eq.(\ref{e:poisson}) 
\begin{equation} \label{e:smoothsource}
S(r,\theta,\varphi) = 4(r^2-2+3z^2)\text{e}^{-r^2-z^2},
\end{equation}
which admits the following solution
\begin{equation} \label{e:smoothsolution}
f(r,\theta,\varphi) = \text{e}^{-r^2-z^2}.
\end{equation}

We compute the error between the analytical and numerical solution by looking at the maximum difference on the collocation points in each domain, which gives an estimation of the corresponding $L^2_w$ norm of the difference. Results are displayed in Fig. 3 and show as expected an exponential decay of the error as a function of the number of radial collocations points. 

One can compare the efficiency of Chebyshev and Jacobi polynomials for the solution (\ref{e:smoothsource})-(\ref{e:smoothsolution}). When \textsc{lorene} uses Chebyshev base in the nucleus, the mapping is $r = Rx$, where $R$ is the radius of the nucleus, and we use the parity of $f_{lm}$ with respect to $l$ to use odd or even Chebyshev bases. Notice this isn't possible in the Jacobi case because Jacobi polynomials with $\alpha \neq \beta$ do not have a definite parity. Results are displayed in Fig. 4. One can see no loss or gain of precision is acquired with this method. In fact, Jacobi polynomials are expected to be more useful in dynamical evolutions to ensure $L^2$ stability (see Sec.~\ref{motivation}).

Finally, we examined a non-smooth solution, namely
\begin{equation}
\begin{array}{l l c l r}
S = 35\sqrt{r}/4  & \text{for} & r \leq R, & S=0 & \text{otherwise}
\\ f=r^{5/2}-7R^{5/2}/2 & \text{for} & r \leq R, & f = -5R^{7/2}/2r & \text{otherwise}
\end{array}
\end{equation}
where $R$ is the outer radius of the shell.

\begin{figure}[t]
\centerline{\includegraphics[width=0.85\columnwidth]{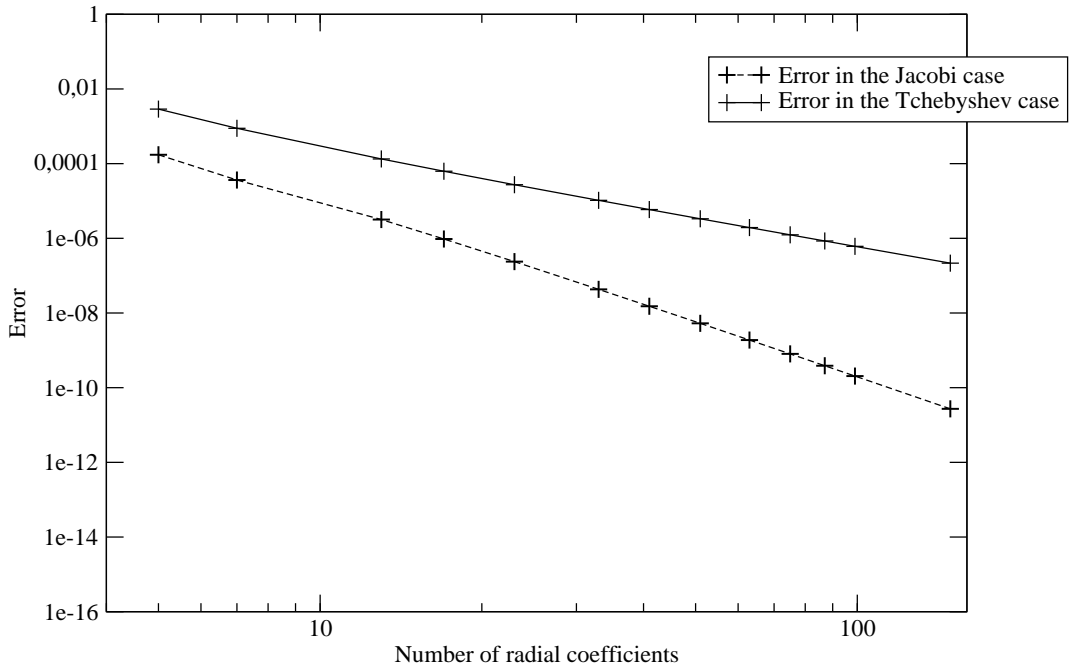}}
\caption[]{\label{f:CV_sqrt} Decay of the errors (maximum difference with theoretical
  solution on collocation points) in the nucleus for Jacobi and Chebyshev methods. Settings are $N_{\theta}=17, N_{\varphi}=16$.}
\end{figure}

Results are displayed in Fig. 5, and one can see as expected a geometric decay of the error as a function of the number of radial coefficients. We found that the rate of convergence in the Chebyshev case is approximately $2.76$, and $4.62$ in the Jacobi case. Indeed, $f \in H^{2}_{w=(1-x^2)^{-1/2}}$, and $f \in H^{3}_{w=(1+x)^2}$.

\section{Concluding remarks}
\label{s:concl}

We have described a new numerical spectral method and showed some applications to Poisson problems in a sphere. In order to build this method, we have first recalled basic results about orthogonal polynomials and generalized quadrature formulas and discrete polynomial tranforms. Besides, we computed various differential and algebraic operators.  Through numerical tests of Poisson equations, we have shown that the method is convergent and accurate. In particular, for smooth solutions, we find as expected a error decaying exponentially with the number of radial coefficients. This method has also been integrated to a 3D multi-domain spectral solver. 
However, the interpolation polynomial error on the Gauss-Lobatto nodes remains to be computed. Studying the nodes of the Gauss Lobatto quadrature formulas shows evidence that their distribution is very similar to the Legendre case, or the $(0,1)$ case, and therefore, that a similar theorem could be proved. Other applications for $(0,2)$ Jacobi polynomials could be foreseen. For instance, hyperbolic equations in a sphere, with applications to fluid dynamics in a star.

\section*{Acknowlegments}
\label{sec:acknowledg}

The authors wish to thank J\'er\^ome Novak for his help on numerical issues, Eric Gourgoulhon and Christine Bernardi for their critical reading of the manuscript. This work was supported by  ANR grant 06-2-134423 \emph{M\'ethodes math\'ematiques pour la relativit\'e g\'en\'erale}.

\appendix
\section{Appendix}

\label{s:annex}
In this section we provide proofs for the major new results listed in this paper.

\subsection{Proposition \ref{p:JGLW}}
\begin{proof}

To simplify the expressions, $J_n$ will be denoted $J_n^{(\alpha,\beta)}$ in this proof, except when the index is explicitly written. The proof resorts to a already known \cite{bernardidaugemaday} preliminary lemma. We shall note $J_{n+1}' = (\lambda_n x + \mu_n)J_n' - \nu_n J_{n-1}'$ the induction formula of the $J_n'$, since they are orthogonal polynomials with respect to the weight $w(x)(1-x^2)$.

Let us first compute $\rho_i$, for $1 \leq i \leq N-1$.

\begin{lem}
\begin{equation} \forall i \in [1,N-1] \quad \rho_i = \frac{\lambda_{n-1}\| J_{n-1}'\|_{L_{w(1-x^2)}^2(\Lambda)}^2} {(1-x_i)^2 J_n''(x_i) J_{n-1}'(x_i) }. \end{equation}
\end{lem}

\begin{proof}[Lemma proof]
Let us fix $i$ in $[1,N-1]$, et define $\Phi(x) = \frac{J_N'(x)}{x-x_i} (1-x^2)$ in $\mathbb{P}_N(\Lambda)$. The exactness of the quadrature formula (\ref{quadformula}) then gives
\begin{equation} (1-x_i^2) \rho_i J_N''(x_i) = \int_{-1}^1 \frac{J_N'(x)}{x-x_i} (1-x^2) w(x) dx .\end{equation}
The computation of the integral may be done through the study of the following quantity
\begin{equation} S_n(\zeta,\eta)= \frac{J_{n+1}'(\zeta)J_n'(\eta) - J_{n+1}'(\eta) J_n'(\zeta)}{\zeta-\eta}. \end{equation}
The induction formula of the $J_n'$ allows us to write
\begin{equation} S_n(\zeta,\eta) = \lambda_n J_n'(\zeta)J_n'(\eta) + \nu_n S_{n-1}(\zeta,\eta). \end{equation}
We here recall that $\nu_n=\lambda_n \| J_n'\|_{L_{w(1-x^2)}(\Lambda)}^2 / \lambda_{n-1} \| J_{n-1}'\|_{L_{w(1-x^2)}(\Lambda)}^2$, according to the induction formula of orthogonal polynomials.
Then
\begin{equation} \frac{S_n(\zeta,\eta)}{\lambda_n \| J_n'\|_{L_{w(1-x^2)}^2(\Lambda)}^2} = \frac{J_n'(\zeta) J_n'(\eta)}{\| J_n'\|_{L_{w(1-x^2)}^2(\Lambda)}^2} + 
\frac{S_{n-1}(\zeta,\eta)}{\lambda_{n-1} \| J_{n-1}' \|_{L_{w(1-x^2)}^2(\Lambda)}^2 }. \end{equation}
Since $S_0(\zeta,\eta)=0$, one can see that
\begin{equation} S_n(\zeta,\eta) = \lambda_n \| J_n'\|_{L_{w(1-x^2)}^2(\Lambda)}^2 \sum_{k=1}^n \frac{J_k'(\zeta) J_k'(\eta)}{\| J_k'\|_{L_{w(1-x^2)}^2(\Lambda)}^2} \end{equation}
This is already known Christoffel-Darboux formula (see \cite{handbook}) for orthogonal polynomials.
Let us take that equality in $n=N-1$, $\zeta = x$, $\eta = x_i$, multiply by $(1-x^2)$, then integrate with respect to $w$. Left member is given by
\begin{equation} \int_{-1}^1 \frac{J_N'(x)J_{N-1}'(x_i) - J_N'(x_i)J_{N-1}'(x)}{x-x_i}(1-x^2)w(x) dx = J_{N-1}'(x_i) \int_{-1}^1 \frac{J_N'(x)}{x-x_i} (1-x^2) w(x) dx. \end{equation}
In the right member, by orthogonality, every term in the sum vanishes, except for  $k=1$, to give
\begin{equation} \lambda_{N-1} \frac{\| J_{N-1}' \|_{L_{w(1-x^2)}^2(\Lambda)}^2   J_1'(x_i)}{\| J_{1}' \|_{L_{w(1-x^2)}^2(\Lambda)}^2} \int_{-1}^1 J_1'(x)(1-x^2) w(x) dx = \lambda_{N-1} \| J_{N-1}' \|_{L_{w(1-x^2)}^2(\Lambda)}^2  \end{equation}
since $J_1'$ is a constant. 
Then, one deduces that
\begin{equation} \int_{-1}^1 \frac{J_N'(x)}{x-x_i} (1-x^2) w(x) dx = \frac{\lambda_{N-1} \| J_{N-1}' \|_{L_{w(1-x^2)}^2(\Lambda)}^2 }{J_{N-1}'(x_i)}. \end{equation}

\end{proof}
\medskip

Let us now compute the various quantities in the lemma.
\begin{itemize}

\item Computation of $\| J_{N-1}' \|_{L_{w(1-x^2)}^2(\Lambda)}^2 $

After integration by parts, and using the differential equation (\ref{eqdifjacobi}), one has

\begin{equation}  \int_{-1}^1 (J_{N-1}')^2 (1-x^2) w(x) dx  =  (N-1)(N+\alpha + \beta) \int_{-1}^1 J_{N-1}^2(x) w(x) dx \end{equation}
i.e,
\begin{equation} \| J_{N-1}' \|_{L_{w(1-x^2)}^2(\Lambda)}^2  = (N-1)(N +\alpha +\beta )\frac{2^{\alpha+\beta+1}}{2N+\alpha+\beta-1}  \frac{\Gamma(N+\alpha)\Gamma(N+\beta)}{\Gamma(N)\Gamma(N+\alpha+\beta)} . \end{equation}

\item Computation of $\lambda_{N-1}$

The identification of highest degree coefficients leads to
\begin{equation} \lambda_{N-1} =\frac{N k_{N}}{ (N-1)k_{N-1}} = \frac{(2N+\alpha+\beta)(2N-1+\alpha+\beta)}{2(N-1)(N+\alpha+\beta)} .\end{equation}

\item Computation of $(1-x_i^2)J_N''(x_i)$

Developing the differential equation (\ref{eqdifjacobi}), and evaluating in $x_i$ results in

\begin{equation} (1-x_i^2)J_N''(x_i) = -N(N+\alpha+\beta+1) J_N(x_i). \end{equation}

\item Computation of $J_{N-1}'(x_i)$

Induction formula of $J_n'$ evaluated in $x_i$ gives 
\begin{equation} J_{N+1}'(x_i) = - \nu_N J_{N-1}'(x_i). \end{equation}

Besides, deriving the induction formula (\ref{inducjacobi}) of the $J_n$, and evaluating in $x_i$ gives  

\begin{multline}  J_{N+1}'(x_i) = \frac{(2N+1+\alpha+\beta)(2N+2+\alpha+\beta)}{2(N+1)(N+1+\alpha+\beta)} J_N(x_i)  
\\
 - \frac{(N+\alpha)(N+\beta)(2N+2+\alpha+\beta)}
{(N+1)(N+1+\alpha+\beta)(2N+\alpha+\beta)} J_{N-1}'(x_i).
\end{multline}

But
\begin{equation} \nu_N = \frac{\lambda_N \| J_N'\|_{L_{w(1-x^2)}(\Lambda)}^2}{\lambda_{N-1} \| J_{N-1}' \|_{L_{w(1-x^2)}(\Lambda)}^2} = \frac{(N+\alpha)(N+\beta)(2N+2+\alpha+\beta)} {N(N+\alpha+\beta)(2N+\alpha+\beta)}. \end{equation}

After substitution, one has

\begin{multline} J_{N-1}'(x_i) \frac{(N+\alpha)(N+\beta)(2N+2+\alpha+\beta)}{2N+\alpha+\beta}
\\
 \left( \frac{1}{(N+1)(N+\alpha+\beta+1)} - \frac{1}{N(N+\alpha+\beta)} \right) 
\\ = \frac{(2N+2+\alpha+\beta)(2N+1+\alpha+\beta)}{2(N+1)(N+1+\alpha+\beta)} J_N(x_i) .
\end{multline}

So

\begin{equation} J_{N-1}'(x_i) = -\frac{1}{2}\frac{N(N+\alpha+\beta)}{(N+\alpha)(N+\beta)}(2N +\alpha +\beta) J_N(x_i). \end{equation}

\end{itemize}

Collecting all above formulas, one has

\begin{equation} \rho_i = \frac{2^{\alpha+\beta+1}}{N(N+\alpha+\beta+1)} \frac{\Gamma(N+1+\alpha)\Gamma(N+1+\beta)}{\Gamma(N+1)\Gamma(N+\alpha+\beta+1)} \frac{1}{(J_N(x_i))^2}. \end{equation}

\bigskip

Let us now compute $\rho_N$. We apply the quadrature formula (\ref{quadformula}) to $\Phi(x) = J_N'(x) (1+x)$. We then obtain
\begin{equation} 2 \rho_N J_N'(1) = \int_{-1}^{1} J_N'(x) (1+x) w(x) dx = \frac{N+\alpha+\beta+1}{2} \int_{-1}^1 J_{N-1}^{(\alpha+1,\beta+1)}(x) (1+x) w(x) dx. \end{equation}
But $(2N+2+\alpha+\beta)J_N^{(\alpha,\beta+1)} = (N+2+\alpha+\beta) J_N^{(\alpha+1,\beta+1)} - (N+\beta+1) J_{N-1}^{(\alpha+1,\beta+1)}$, according to (\ref{linkjacobi}). Multiplying this by $(1+x)w(x)$ and integrating leads to
\begin{equation} (N+\alpha+\beta+2) \int_{-1}^1 J_N^{(\alpha+1,\beta+1)}(x)(1+x)w(x) dx = (N+\beta+1) \int_{-1}^1 J_{N-1}^{(\alpha+1,\beta+1)} (1+x)w(x) dx. \end{equation}
Hence 
\begin{equation} \int_{-1}^1 J_N^{(\alpha+1,\beta+1)}(x)(1+x)w(x) dx = \frac{\Gamma(N+\beta+2)\Gamma(\alpha+\beta+3)}{\Gamma(\beta+2)\Gamma(N+\alpha+\beta+3)} \int_{-1}^1 J_0^{(\alpha+1,\beta+1)}(x)(1+x)w(x)dx.
\end{equation}
Moreover,
\begin{equation} \int_{-1}^1 (1+x)w(x) dx = \| J_0^{(\alpha,\beta+1)} \|_{L_{(1+x)w}^2(\Lambda)}^2 = 2^{\alpha+\beta+2} \frac{\Gamma(\alpha+1)\Gamma(\beta+2)}{\Gamma(\alpha+\beta+3)}. \end{equation}
Finally,
\begin{equation} J_N'(1) = \frac{1}{2}(N+\alpha+\beta+1) \frac{\Gamma(N+1+\alpha)}{\Gamma(2+\alpha)\Gamma(N)}. \end{equation}
Then,
\begin{equation} \rho_N = 2^{\alpha+\beta+1} \frac{\Gamma(\alpha+2)\Gamma(\alpha+1)\Gamma(N)\Gamma(N+1+\beta)} {\Gamma(N+1+\alpha)\Gamma(N+2+\alpha+\beta)}. \end{equation}
Since $J_N(1)=\frac{\Gamma(N+1+\alpha)}{\Gamma(1+\alpha)\Gamma(N+1)}$, one has
\begin{equation} \rho_N = (\alpha+1) \frac{2^{\alpha+\beta+1}}{N(N+\alpha+\beta+1)} \frac{\Gamma(N+1+\alpha)\Gamma(N+1+\beta)}{\Gamma(N+1)\Gamma(N+1+\alpha+\beta)} \frac{1}{(J_N^{(\alpha,\beta)}(x_n))^2}. \end{equation}

The expression of $\rho_0$ is obtained by exchanging $\alpha$ and $\beta$.

\end{proof}

\subsection{Proposition \ref{p:matrixnodes}}
\begin{proof}
First, let us recall that $\frac{d}{dx} \left( J_n^{(\alpha,\beta)} \right) = \frac{1}{2}(n+\alpha+\beta+1) J_{n-1}^{(\alpha+1,\beta+1)}$, then we will look for the zeros of $J_{N-1}^{(\alpha+1,\beta+1)}$. But the $J_n^{(\alpha+1,\beta+1)}$ satisfy the following induction formula
\begin{multline} 2 (n+1) (n+\alpha+\beta+3) (2n+\alpha+\beta+2) J_{n+1}^{(\alpha+1,\beta+1)} = \\
		\left[ (2n+\alpha+\beta+3) (\alpha-\beta) (\alpha+\beta+2) + x \Gamma(2n+\alpha+\beta+5)/\Gamma(2n+\alpha+\beta+2) \right] J_n^{(\alpha+1,\beta+1)}  \\
		- 2 (n+\alpha+1) (n+\beta+1) (2n+\alpha+\beta+4) J_{n-1}^{(\alpha+1,\beta+1)} 
\end{multline}
Let us then define
\begin{equation} J_n^* = \frac{J_n^{(\alpha+1,\beta+1)}}{\|J_n^{(\alpha+1,\beta+1)}\|_{L_{(1-x^2)w}^2(\Lambda)}} = \sqrt{\frac{(2n+\alpha+\beta+3) \Gamma(n+1) \Gamma(n+\alpha+\beta+3)} {2^{\alpha+\beta+3} \Gamma(n+\alpha+2) \Gamma(n+\beta+3)} }J_n^{(\alpha+1,\beta+1)}  \end{equation}
Then, the reccurence formula of the $J_n^*$ can be written as such
\begin{multline}  
x J_{n-1}^* =-\frac{(\alpha-\beta)(\alpha+\beta+2)}{(2n+\alpha+\beta)(2n+\alpha+\beta+2)} J_{n-1}^* \\
	+ \frac{2}{2n+\alpha+\beta+2} \sqrt{ \frac{n(n+\alpha+1)(n+\beta+1)(n+\alpha+\beta+2)} {(2n+\alpha+\beta+1)(2n+\alpha+\beta+3)}} J_n^* \\
	+ \frac{2}{2n+\alpha+\beta} \sqrt{ \frac{(n-1)(n+\alpha)(n+\beta)(n+\alpha+\beta+1)} {(2n+\alpha+\beta-1)(2n+\alpha+\beta+1)}} J_{n-2}^* 
\end{multline}
or similarly
\begin{equation} x J_{n-1}^* = \delta_n J_{n-1}^* + \gamma_n J_n^* + \gamma_{n-1} J_{n-2}^* \end{equation}
which can be written in a matrix form
\begin{equation} x \left( \begin{array}{c}  J_0^*(x) \\ J_1^*(x) \\ \vdots \\ J_{N-3}^*(x) \\ J_{N-2}^*(x) \end{array} \right) = \left( \begin{array}{c c c c c}
\delta_1 & \gamma_1 & 0        & 0            & 0            \\
\gamma_1 & \delta_2 & \ddots   & 0            & 0            \\
0        & \ddots   & \ddots   & \ddots       & 0            \\
0        & 0        & \ddots   & \delta_{N-2} & \gamma_{N-2} \\
0        & 0        & 0        & \gamma_{N-2} & \delta_{N-1}
 \end{array} \right) \left( \begin{array}{c}  J_0^*(x) \\ J_1^*(x) \\ \vdots \\ J_{N-3}^*(x) \\ J_{N-2}^*(x) \end{array} \right) + \gamma_{N-1} \left( \begin{array}{c}  0 \\ 0 \\ \vdots \\ 0 \\ J_{N-1}^*(x) \end{array} \right) \end{equation}
Otherly said, the $x_j$, $1 \leq j \leq N-1$, zeros of $J_{N-1}^*$, are the eigenvalues of the (\ref{e:matrixnodes}) matrix.

\end{proof}

\subsection{Proposition \ref{p:djt}}
\begin{proof}

One can write, for a fixed $m \leq N$, 
\begin{equation} \sum_{i=0}^N \rho_i J_m(x_i) i_Nf(x_i) = \sum_{k=0}^N \tilde{f}_k \left\{ \sum_{i=0}^N  \rho_i J_k(x_i) J_m(x_i) \right\}. \end{equation}

For $m \leq N-1$, since $k+m \leq 2N-1$, only $ \tilde{f}_m \|J_m\|_{L_w^2(\Lambda)}^2$ is left in the right-hand side. Then $ \tilde{f}_m = \frac{1}{\|J_m\|_{L_w^2(\Lambda)}^2} \sum_{i=0}^N \rho_i J_m(x_i)f(x_i)$. Eq.(\ref{e:jacobinorm}) and the expressions of $\rho_i$ from proposition \ref{p:JGLW} lead to the desired result.

For $m = N$, the exactness of the quadrature formula allows us to suppress the $N-1$ first terms of the sum, but the last one is not equal to the norm of $J_N$. Precisely, the right side member is equal to $\tilde{f}_N \sum_{i=0}^N \rho_i (J_N(x_i))^2$. Using the expressions of the weights, one finds

\begin{equation} \tilde{f}_N \frac{2^{\alpha+\beta+1}}{N(N+\alpha+\beta+1)} \frac{\Gamma(N+1+\alpha)\Gamma(N+1+\beta)}{\Gamma(N+1)\Gamma(N+1+\alpha+\beta)} \left[ (\beta+1) +(N-1) + (\alpha+1) \right]. \end{equation}
One deduces 
\begin{equation} \tilde{f}_N \frac{2^{\alpha+\beta+1}}{N} \frac{\Gamma(N+1+\alpha)\Gamma(N+1+\beta)}{\Gamma(N+1)\Gamma(N+1+\alpha+\beta)} = \sum_{i=0}^N \rho_i J_N(x_i)f(x_i), \end{equation}
which is the desired result.
\end{proof}

\subsection{Theorem \ref{t:invineq}}
\begin{proof}

First let us compute $\| J_n' \|_{L_w^2(\Lambda)}^2$ thanks to the Gauss-Lobatto quadrature with $n+1$ points ($J_n'^2$ has indeed degree $2n-2$).

\begin{multline} \int_{-1}^1 (J_n'(x))^2 w(x) dx =
\rho_0 (J_n'(-1))^2 + \rho_n (J_n'(1))^2 
\\
= \frac{(n+\alpha+\beta+1)^2}{4} \left( \rho_0 \left( \frac{\Gamma(n+\beta+1)}{\Gamma(n)\Gamma(\beta+2)} \right)^2 + \rho_N \left( \frac{\Gamma(n+\alpha+1)}{\Gamma(n)\Gamma(\alpha+2)} \right)^2 \right) 
\\
= 2^{\alpha+\beta-1}n(n+\alpha+\beta+1)\left(\frac{1}{\alpha+1} + \frac{1}{\beta+1}\right) \frac{\Gamma(n+1+\alpha)\Gamma(n+1+\beta)}{\Gamma(n)\Gamma(n+1+\alpha+\beta)} .
\end{multline} 
Then, for $\varphi_N = \sum_{n=0}^N \varphi^n J_n$, one has, thanks to Cauchy-Schwarz inequality
\begin{equation}
|\varphi_N|_{H^1_w(\Lambda)} \leq \sum_{n=0}^{N} |\varphi^n| |J_n|_{H^1_w(\Lambda)} 
\leq \left( \sum_{n=0}^N |\varphi^n|^2 \|J_n\|_{L^2_w(\Lambda)}^2 \right)^{1/2} \left( \sum_{n=0}^N \frac{|J_n|^2_{H^1_w(\Lambda)}}{\|J_n\|^2_{L_w^2(\Lambda)}} \right)^{1/2}.
\end{equation}
But, according to the previous computation,

\begin{equation} \frac{|J_n|^2_{H^1_w(\Lambda)}}{\|J_n\|^2_{L_w^2(\Lambda)}} = \frac{1}{4} \left( \frac{1}{\alpha+1}+\frac{1}{\beta+1} \right) n(n+\alpha+\beta+1)(2n+\alpha+\beta+1). \end{equation}
One can deduce that \begin{multline}  |\varphi_N|_{H^1_w (\Lambda)} \leq  \| \varphi \|_{L^2_w(\Lambda)} \left[ N^3 \frac{1}{4}\left(\frac{1}{\alpha+1}+\frac{1}{\beta+1}\right) \sum_{n=0}^N \frac{n}{N}\frac{n+\alpha+\beta+1}{N}\frac{2n+\alpha+\beta +1}{N} \right]^{1/2} 
\\
\leq 
\| \varphi \|_{L^2_w(\Lambda)} \left[ N^4 \frac{1}{4}\left( \frac{1}{\alpha+1}+\frac{1}{\beta+1}\right)  \text{sup}_{0\leq n \leq N } \left(\frac{n}{N}\frac{n+\alpha+\beta+1}{N}\frac{2n+\alpha+\beta+1}{N}\right)   \right]^{1/2}.
\end{multline}
Therefore, we have \begin{equation} |\varphi_N|_{H^1_w(\Lambda)} \leq c N^2 \| \varphi_N \|_{L^2_w(\Lambda)}. \end{equation}

Now, we turn to the optimality of this inequality. We will compute $\| J_N'' \|_{L_w^2(\Lambda)}$ thanks to the quadrature formula with $N+1$ points.

\begin{equation}
\| J_N'' \|_{L_w^2(\Lambda)}^2 \geq \rho_0 (J_N''(-1))^2 + \rho_N (J_N''(1))^2
\end{equation}
because the weights $\rho_i$ are non-negative. Using (\ref{derjacobi}) and (\ref{e:normalisationjacobi}), one can find that
\begin{multline}
\| J_N'' \|_{L_w^2(\Lambda)}^2 \geq 2^{\alpha+\beta-3}N(N-1)^2(N+\alpha+\beta+1)(N+\alpha+\beta+2)^2 \\ \frac{\Gamma(N+1+\alpha)\Gamma(N+1+\beta)}{\Gamma(N+1)\Gamma(N+1²+\alpha+\beta)} \left( \frac{1}{(\beta+1)(\beta+2)^2} + \frac{1}{(\alpha+1)(\alpha+2)^2} \right).
\end{multline}
Thus,
\begin{equation}
\frac{\| J_N'' \|_{L_w^2(\Lambda)}^2}{\| J_N' \|_{L_w^2(\Lambda)}^2}  \geq \frac{1}{4} (N-1)^2(N+\alpha+\beta+2)^2 \left( \frac{\alpha+1}{(\beta+2)^2} + \frac{\beta+1}{(\alpha+2)^2} \right) \frac{1}{\alpha+\beta+2}.
\end{equation}
Therefore,
\begin{equation} 
\frac{| J_N' |_{H_w^1(\Lambda)}^2}{\| J_N' \|_{L_w^2(\Lambda)}^2} \geq c' N^4
\end{equation}

\end{proof}

\subsection{Proposition \ref{p:opmatrices}}
\begin{proof}

\begin{itemize}

\item Derivation :

We recall that \begin{equation} \begin{array}{l} 
(2n+4)J_n^{(0,3)} = (n+4)J_n^{(1,3)} - (n+3)J_{n-1}^{(1,3)} \\
(2n+3)J_n^{(0,2)} = (n+3)J_n^{(0,3)} + n J_{n-1}^{(0,3)} 
\end{array} \end{equation}
From these equalities, we obtain \begin{equation} \begin{array}{l}
(n+4)J_n^{(1,3)} = 2 \sum_{k=0}^n (k+2) J_k^{(0,3)} \\
(n+3)(n+2)(n+1)J_n^{(0,3)} = (-1)^n \sum_{k=0}^n (-1)^k (2k+3)(k+1)(k+2)J_k^{(0,2)} 
\end{array}, \end{equation}
which allows us to write \begin{equation} (n+4)J_n^{(1,3)} = 2\sum_{j=0}^n (2j+3)(j+1)(j+2)(-1)^j \left\{  \sum_{k=j}^n \frac{(-1)^k}{(k+3)(k+1)} \right\} J_j^{(0,2)}. \end{equation}
The term in between brackets can be computed using a simple element decomposition and leads to \begin{equation} 2 \sum_{k=j}^n \frac{(-1)^k}{(k+3)(k+1)} = \frac{(-1)^n}{(n+2)(n+3)} + \frac{(-1)^j}{(j+1)(j+2)}. \end{equation}
One can deduce that \begin{equation} (n+4)J_n^{(1,3)} = \sum_{j=0}^n (2j+3) \left\{ 1 + (-1)^{j-n}\frac{(j+1)(j+2)}{(n+2)(n+3)} \right\} J_j^{(0,2)}, \end{equation}
Which gives the desired result, thanks to $\frac{d}{dx}J_n^{(0,2)} = \frac{1}{2}(n+3)J_{n-1}^{(1,3)}$.

\bigskip

\item Integration :

Let us recall that
\begin{equation} \begin{array}{l}
(2n+1)J_n^{(-1,1)} = (n+1)J_n^{(0,1)} - (n+1)J_n^{(0,1)} \\
2(n+1)J_n^{(0,1)} = (n+2) J_n^{(0,2)} + n J_{n-1}^{(0,2)}
\end{array} \end{equation}
Then, \begin{multline} 
\frac{1}{2}(k+1) J_{k-1}^{(0,2)} 
= \frac{d}{dx} J_k^{(-1,1)} 
= \frac{k+1}{2k+1} \frac{d}{dx} (J_k^{(0,1)} - J_{k-1}^{(0,1)}) \\ 
= \frac{k+1}{2k+1} \frac{d}{dx} \left( \frac{(k+2) J_k^{(0,2)}+k J_{k-1}^{(0,2)}}{2(k+1)} -  \frac{(k+1) J_{k-1}^{(0,2)}+(k-1) J_{k-2}^{(0,2)}}{2k} \right). 
\end{multline}
Or \begin{equation} J_{k-1}^{(0,2)} = \frac{d}{dx} \left( \frac{k+2}{(k+1)(2k+1)} J_k^{(0,2)} - \frac{1}{k(k+1)} J_{k-1}^{(0,2)} - \frac{k-1}{k(2k+1)}J_{k-2}^{(0,2)} \right). \end{equation} By taking $n=k-1$, we obtain the primitive given by the propostion \ref{p:opmatrices}, besides, since $J_n^{(0,2)}(1)=1$, this primitive vanishes at $1$. 

\item Division by $1+x$ :

We will use the following formula :
\begin{equation} (n+2)(1+x)J_n^{(0,3)} = (n+3)J_n^{(0,2)} + (n+1)J_{n+1}^{(0,2)}, \end{equation}
from which \begin{equation} \frac{J_{n+1}^{(0,2)}}{1+x} = \sum_{k=0}^n (-1)^{n-k} \frac{(n+3)(n+2)}{(k+3)(k+2)(k+1)} (k+2)J_k^{(0,3)} + \frac{J_{n+1}^{(0,2)}(-1)}{1+x}. \end{equation}
But we saw previously that \begin{equation} (n+3)(n+2)(n+1)J_n^{(0,3)} = (-1)^n \sum_{k=0}^n (-1)^k (2k+3)(k+1)(k+2)J_k^{(0,2)}. \end{equation} 
So \begin{multline} \frac{J_{n+1}^{(0,2)}}{1+x} = (n+3)(n+2)\sum_{j=0}^n (-1)^{n-j}(2j+3)(j+2)(j+1) \\ \left\{ \sum_{k=j}^n \frac{1}{(k+1)^2(k+2)(k+3)^2} \right\} J_j^{(0,2)}  + \frac{J_{n+1}^{(0,2)}(-1)}{1+x}. \end{multline}
The term in between brackets can be computed using the following simple element decomposition:
\begin{equation} \frac{4}{(k+1)^2(k+2)(k+3)^2} = -\frac{2}{k+1} + \frac{1}{(k+1)^2} + \frac{1}{k+2} - \frac{2}{k+3} - \frac{1}{(k+3)^2}. \end{equation}
The sum presents lots of cancellations and is equal to
\begin{equation} 4\sum_{k=j}^n \frac{1}{(k+1)^2(k+2)(k+3)^2} = \frac{1}{(j+1)^2(j+2)^2} - \frac{1}{(n+2)^2(n+3)^2}. \end{equation}
Finally, \begin{equation} \frac{J_{n+1}^{(0,2)}}{1+x} = \sum_{j=0}^n (-1)^{n-j}\frac{2j+3}{4} \left\{ \frac{(n+3)(n+2)}{(j+1)(j+2)} - \frac{(j+2)(j+1)}{(n+2)(n+3)} \right\} J_j^{(0,2)} + \frac{J_{n+1}^{(0,2)}(-1)}{1+x}. \end{equation}
Notice that this proof can also be achieved using a Christoffel-Darboux formula.

\item Multiplication by $1+x$ :

We recall that \begin{equation} \begin{array}{l} (n+3/2)(1+x)J_n^{(0,2)} = (n+2)J_n^{(0,1)} + (n+1)J_{n+1}^{(0,1)} \\
2(n+1)J_n^{(0,1)} = (n+2) J_n^{(0,2)} + n J_{n-1}^{(0,2)} \end{array} \end{equation}
Then, \begin{multline} (n+3/2)(1+x)J_n^{(0,2)} = (n+2)\left( \frac{(n+2)J_n^{(0,2)} + nJ_{n-1}^{(0,2)}}{2(n+1)} \right) \\+ (n+1)\left( \frac{(n+3)J_{n+1}^{(0,2)} + (n+1)J_n^{(0,2)}}{2(n+2)} \right). \end{multline}
Or \begin{equation} (1+x)J_n^{(0,2)} = \frac{(n+1)(n+3)}{(n+2)(2n+3)}J_{n+1}^{(0,2)} + \frac{n^2+3n+3}{(n+1)(n+2)}J_n^{(0,2)} + \frac{n(n+2)}{(n+1)(2n+3)}J_{n-1}^{(0,2)}. \end{equation}

\end{itemize}
\end{proof}

\end{document}